\documentclass[a4paper,10pt]{article}
\usepackage{amsmath,amssymb, amsfonts, textcomp, eufrak}
\usepackage[affil-it]{authblk}
\newtheorem{theorem}{Theorem}[section]
\newtheorem{definition}[theorem]{Definition}
\newtheorem{lem}[theorem]{Lemma}
\newtheorem{prop}[theorem]{Proposition}
\newtheorem{rem}[theorem]{Remark}

\numberwithin{equation}{section}
\usepackage{graphicx}
\numberwithin{equation}{section}
\usepackage{graphicx}
\usepackage{float} 
\usepackage{hyperref}
\usepackage{fancyhdr}
\usepackage{makeidx}
\makeindex
\newcommand{\R}{\mathbb{R}}

\newcommand{\N}{\mathbb{N}}

\newcommand{\Rp}{\mathbb{R}^+}
\newcommand{\F}{\mathcal{F}}

\newcommand{\Bb}{\mathcal{B}}

\newcommand{\K}{(\F)_{t\ge 0}}

\newcommand{\Z}{\mathbb{Z}}

\newcommand{\T}{\mathbb{T}}

\newcommand{\p}{\mathbb{P}}
\newcommand{\Pp}{\mathcal{P}_{\p}(\Omega\times X)}
\newcommand{\bn}{\begin{definition}}
\newcommand{\en}{\end{definition}} 
\newcommand{\bt}{\begin{theorem}}                
\newcommand{\et}{\end{theorem}}
 \newcommand{\bnm}{\begin{enumerate}}              
\newcommand{\enm}{\end{enumerate}}
\newcommand{\br}{\begin{rem}} 
\newcommand{\er}{\end{rem}}

\newcommand{\om}{\omega}

\newcommand{\Om}{\Omega}

\newcommand{\btm}{\begin{itemize}}
 \newcommand{\etm }{\end{itemize}}

\newcommand{\pds}{(\Om\times X, \F\otimes\Bb)}
\newcommand{\pro}{\mathcal{P}_{\p}(\Om\times X)}
\newcommand{\Ip}{\mathcal{I}_{\p}(\varphi)}
\newcommand{\E}{\mathbb{E}}

\newcommand{\Rd}{\R^d}

\newcommand{\hq}{\hat{\theta}}
\newcommand{\Hha}{\hat{H}}

\newcommand{\Ht}{\hat{t}}
\newcommand{\SM}{\mathcal{M}_{\p}}
\newcommand{\Lpr}{\lambda^{\prime}}

\newcommand{\Sc}{\mathbb{S}^1}	

\begin{document}
\markboth{K. Uda}{Existence of random invariant periodic curves}
\title{Existence of random  invariant periodic curves via random semiuniform ergodic theorem }
\author{Kenneth Uda}

\affil{Department of Mathematical Sciences Loughborough University, LE11 3TU, UK \\K.O.Uda@lboro.ac.uk}
\date{}
\providecommand{\keywords}[1]{\textbf{\textit{Keywords:}} #1}
\maketitle
\begin{abstract}\rm We employ an extension of ergodic theory to the random setting to investigate the existence of random periodic solutions of random dynamical systems.  Given that a random dynamical system has a dissipative structure, we proved that a random invariant compact set can be expressed as a union of finite of number of random  periodic curves. The idea in this paper is closely related to the work  recently considered by Zhao and Zheng \cite{Zhao}.
\end{abstract}


\section{Introduction}
Long time behaviour of mathematical models are motivated by two main issues with theoretical and practical consequences. The first is to understand where the orbits (collection of solutions) accumulate in the long run. The second and equally important one is to ascertain whether the limiting behaviour is still essentially the same after small changes to the evolution rule.  Intuitively, the limiting behaviour of dynamical system is captured by the concept of stationary and periodic solutions. For a dynamical system $\varphi:\T\times X\rightarrow X,$ over $t\in\T,$ a stationary solution is a point $x\in X,$ such that \begin{equation}
   \varphi(t,x) = x, \quad \text{for all $t$}\in\T.
   \end{equation}
And a periodic solution is a periodic function $u:\T\rightarrow X$ with period $\tau\neq 0,$ such that \begin{equation}
u(t+\tau) = u(t), \quad \text{and}\quad \varphi(t,u(s))= u(t+s), \quad \text{for all $t,s$}\in\T.
\end{equation} 
To understand and give the existence of such solutions have attracted vast interest in theory and applications. Periodic solutions have been crucial in the qualitative theory of dynamical systems and its systematic consideration was initiated by Poincar\'e in his work \cite{Ponc}. Periodic solutions have been studied for many fascinating physical problems, examples: van der Pol equations \cite{Van}, Lienard equations \cite{Lien}, etc. However, once noise is added, the dynamics start to depend on both time and the noise path, so the above defintion of steady state (stationary and periodic) solutions may not exist for randomly perturbed systems.  

The long time  behaviour of systems become more interesting and difficult when we include noise in the systems, which in numerous applications are unavoidable. The long time behaviour of random dynamical systems is relatively new area of mathematical research and has seen a tremendous progress in the last three decades. As in the dynamical systems setting, random stationary solutions are central in the long time behaviour of random dynamical systems. For a random dynamical system $\varphi: \T\times\Om\times X\rightarrow X,$ over a metric dynamical system $(\Om, \F,\p,(\theta_t)_{t\in\T}),$ a random stationary solution is an $\F$-measurable random variable $Y:\Om\rightarrow X$ such that \begin{equation}
\varphi(t,\om,Y(\om)) = Y(\theta_t\om), \quad \text{for all $t$}\in\T, \quad \p-\text{almost surely.}
\end{equation}
The notion of random stationary solution of a random dynamical system is a natural extension of fixed point solution of the deterministic system. It is a " one force , one solution" setting that describes the pathwise invariance of the stationary solution over time along the dynamic $\theta$ and the pathwise limt of the random dynamical system. The study of stability once the stationary solution is known were motivated by the poineer work of Has\'minskii in 1969 (translated to English in 1980) \cite{Hasm},  other contributions were made by Kushner \cite{Kushn}, Pinksy \cite{Pinks}, Mao \cite{Mao} and others.  Finding such solutions for stochastic evolution equations is one of the basic problems in stochastic analysis and has been studied recently by Schmalfuss \cite{Schmalfuss, Schm}, Caraballo, Kloeden and Schmalfuss \cite{Caraballo}, E, Khanin, Mazel and Sinai \cite{Sinai}, Sinai \cite{Sina1, Sina2},  Zhang and Zhao \cite{Huaiz}. Stable and unstable manifold was also recently considered by Duan, Lu and Schmalfuss \cite{Duan}, Mohammed, Zhang and Zhao \cite{Saleh} and others. 

Analogous to the periodic solutions of dynamical system, the notion of random periodic solutions play similar role to random dynamical systems. In physical world around us (e.g. biology, chemical reactions, climatic dynamics, finance, etc), we encounter many phenomena which repeat after certain interval of time.  Due to the unavoidable random influences, the phenomena may be best described by random periodic solutions rather than periodic solutions. For example, the maximum daily temperature in any particular region is a random process, however, it certainly has periodic nature driven by divine clock due to the rotation of the earth around the sun.
There have been few attempts in physics to study random perturbation of limit cycle for some time (\cite{Kurr}, \cite{Kampen}, \cite{Jeffrey}). One of the challanges that hinders real progress was lack of a rigorous mathematical definition of random periodic solution and appropriate mathematical tools. For a random path with some periodic property, it is unclear what a reasonable mathematical relation between the random position $\psi(t,\om)$ at time $t$ and $\psi(t+\tau,\om)$ at time $t+\tau$ after a period $\tau$ should be. However, as $\psi(t,\om)$ is a true path, so it is not necessarily true that $\psi(t,\om) = \psi(t+\tau,\om).$ To require that $\psi(t+\tau,\om)$ is in the neighbourhood of $\psi(t,\om)$ by considering a small noise perturbation was worthwhile attempt. However, this approach does not apply to many stochastic differential equations and also lack rigour, and the scope of application is limited.  Recently, in the works of Zhao and Zheng \cite{Zhao}, Feng, Zhao and Zhou \cite{Zhou}, Feng and Zhao \cite{Feng}, it has been observed that for fixed $t,$ $(\psi(t+k\tau,\om))_{k\in\Z}$ should be a random stationary solution of the discrete RDS $\varphi(k\tau,\om)$. This then led to the rigorous definition of random periodicity $\psi(t+\tau,\om)= \psi(t,\theta_{\tau}\om).$  For a random dynamical system $\varphi$ over a metric dynamical systems $(\Om,\F,\p,(\theta_t)_{t\in\T}),$ a \textbf{random periodic solution}  is an $\F$-measurable function $\psi: \T\times\Om\rightarrow X,$ of period $\tau$ such that \begin{equation}\label{intreq}
\psi(t+\tau,\om)= \psi(t,\om) \quad \text{and}\quad \varphi(t,\om, \psi(s,\om)) = \psi(t+s, \theta_t\om), \quad \text{for all $t,s$}\in\T.
\end{equation}
The study of random periodic solution is more fascinating and difficult than the deterministic periodic solution. The extra essential difficulty is from the fact that trajectory (solution path) of the random dynamical systems starting at a point on the periodic curve does not follow the periodic curve, but moves from one periodic curve to another one corresponding to different $\om.$ If one considers family of tajectories starting from different points on the closed curve $\psi(.,\om),$ then the whole family of trajectories at time $t\in\T$ will lie on a closed curve corresponding to $\theta_t\om.$  There are few numerical evidence in the literature sugessting the existence of random periodic curves (examples; Stochstic Van der pol oscillator \cite{SDER}, Stochastic Goodwin-Lotka- Volterra model \cite{Goodwin}, Stochastic speculative financial model \cite{Carl}, Stochastic climate dynamics \cite{Mick}).

Our interest in this paper is to investigate the existence of random periodic solutions using geometric approach. It is subtle if not impossible to represent the solutions of stochastic systems in a closed form. So geometric approach becomes a natural or feasible point of view in  the investigation of some important features of random dynamical systems. Our results are based on some geometric invariant structures of RDS and it is mainly based on the extension of ergodic theorem to the random setting, in particular random semiuniform ergodic theorem \cite{Stark13}. The aim is to present results directly accessible to the theory of random bifurcation (dynamical bifurcation), which has numerous problems yet to be explored. In random D-bifurcation, random invariant structures such as  random invariant sets, random invariant measures and Lyapunov exponents are mostly the major objects of investigation. The basic assumptions in this work, are on these invariant structures and these assumptions are quite natural for important systems. 

Recently, Zhao and Zheng \cite{Zhao} used geometric approach to prove the existence of random periodic solutions on the cylinder $\Sc\times\R^d,$ to the best of our knowledge, this is the first result on the geometric approach to the existence random periodic solutions. It is assumed in \cite{Zhao} that there is a random invariant compact set which consists of Lipschitz continuous graph and it is proved that the graph is random periodic. In fact we did not assume the that the random invariant consists of Lipschitz graph, but we employed random semiuniform ergodic theorem to prove that the random invariant compact set consists of random periodic graph. In fact, we gave some conditions on the bound of top Lyapunov exponent of the RDS in the neighbourhood of the random invariant compact set to achieve this result. This random periodic graph is  continuous (Theorem \ref{cong}) and we are not sure if it is Lipschitz continuous, this would be investigated in the future work. As mentioned in the previous paragraph, the conditions here are consistent with that of random D-bifurcation (\cite{Arnold}, \cite{Carl}, \cite{Hoppe}, \cite{Hoppe1} and \cite{Wangf}) and we believe that further systematic study of our result would be useful in the Hopf bifurcation theory of random dynamical systems.

Finally, let us outline the structure of this paper, in section 2, we present random ergodic theory relevant in this work most of which is from the work of J\"ager and Keller \cite{Stark13}, Sturman and Stark \cite{Stark} and the Arnold's book \cite{Arnold}. In section 3, we present the existence of random periodic solutions on the cylinder $\Sc\times \R^d$ which is the main results of this work. In particular, the result would be applied to deduce the existence of random periodic solutions of the following time periodic stochastic differential equation
\begin{equation}\label{1c1}
dX(t) = F(t,X(t))dt+ G(t,X(t))dW_t.
\end{equation}
For this, we put the stochastic differential equation (\ref{1c1}) in autonomous form on $\R^{d+1}$ namely;
\begin{equation}\label{aut}
\begin{cases}
dX(r)= F(r(t),X(t))dr+ G(r(t),X(t))dW_t,\\
dr(t) = dt,\\
X(0)= x, \quad  r(0)= t_0.
\end{cases}
\end{equation}
 Due to the $\tau$ periodicity of the functions $F$ and $G,$ we can identify each point $(t,x)\in\R^{d+1}$ by $(t+\tau, x),$ hence we have time homogeneous equation on $I_{\tau}\times\R^d,$ where $I_{\tau}$ is the closed interval with end point identified, i.e., $I_{\tau}$ is a circle $\Sc.$ Technically, here we have $s+t\mod\tau,$ i.e., $k\in\N$
$$
s+t-k\tau= \tau(\frac{s+t}{\tau}-k)= \tau (\frac{s+t}{\tau}\mod 1).
$$

\section{Random semiuniform ergodic theorem}
Suppose $\theta$ is an ergodic flow of measure preserving transformations of $(\Om,\F,\p),$ the family of maps $\{\varphi(t,\om)\}$ is called \textbf{random dynamical system} if 
\btm
\item $(t,\om,x)\mapsto \varphi(t,\om,x)$ is measurable,
\item $\varphi(0,\om) = \text{id}$ for all $\om\in\Om$ and 
\item $\varphi(t+s,\om)= \varphi(t,\theta_s\om)\circ\varphi(s,\om)\quad \text{for all}\quad s,t\in\T,$ $\p$-almost surely.
\etm
The notion of random dynamical system comprises  stochastic (partial) differential equations, random differential equations and products of random diffeomorphisms.

Let $X$ be a Polish space, an invariant meausre for a random dynamical system is a probability measure on the product space $\Om\times X$ with marginal $\p$ on $\Om,$ which is invariant under the skew product dynamical system $\Theta:=(\theta,\varphi)$ on $\Om\times X.$ A random probability measure $\mu$ on $\Om\times X$ is uniquely characterized by its factorization $\om\mapsto\mu_{\om},$ where
$$\mu(A) = \int_{\Om\times X}\mathbb{I}_{A}(\om,x)d\mu_{\om}(x)d\p(\om), \quad A\in\F\otimes\Bb,$$ 
also written as
 $$d\mu(x,\om) = d\mu_{\om}(x)d\p(\om).$$
\btm
\item$\pro:= \{\mu \quad \text{probability measure on}\quad \pds,\quad \text{with}\quad \pi_{\Om}\mu = \p \},$
\item$\Ip:= \{\mu\in\pro: \mu \quad \text{is}\quad \varphi\text{-invariant}\}.$ 
\item
$\SM(K):= \big\{\mu\in\Ip: \mu(K)=1\big\}$- the set of all $\varphi$-invariant probability measures supported on the random compact set $K.$
\etm

 Assume that $\Ip$ is nonvoid (see Theorem \ref{Kry} in the appendix), let $\mu\in \Ip$ and let $(\Phi_n)_{n\in\N}$ be a sequence of subadditive functions with respect to $\Theta.$ We define 
\begin{equation}\label{invm}
\mu(\Phi_n):= \int_{\Om\times X}\Phi_n(\om,x)d\mu_{\om}(x)d\p(\om),
\end{equation} 
and subadditivity of $(\Phi_n)_{n\in\N}$ with respect to the skew product $\Theta,$ implies that 
\begin{equation}\label{sub1}
\mu(\Phi_{n+m})\leq \mu(\Phi_n)+\mu(\Phi_m)
\end{equation} 
provided both sides are well defined, such that the subadditivity limit 
\begin{equation}\label{sub2}
\lambda(\mu,\varphi):=\lim_{n\rightarrow\infty}\frac{1}{n}\mu(\Phi_n)
\end{equation}
exists. 
 
 Given a random continuous function $\Phi:\Om\times X\rightarrow \R$ and a random compact set $K,$ we define \begin{equation}\label{max}
 \Phi^{K}(\om):= \max\{\Phi(\om,x): x\in K(\om)\}.
 \end{equation}
 \begin{lem}[\cite{Stark13}]\label{maxly}\rm The function
 $\Phi^K:\Om\rightarrow\R$ is measurable and there is a measurable function $b:\Om\rightarrow X$ such that $b(\om)\in K(\om)$ and $\Phi^K(\om)= \Phi(\om,b(\om))$ for all $\om\in\Om.$ 
 \end{lem}
 
 \begin{lem}[\cite{Stark13}, \cite{Cao}]\label{26}\rm
 
 Let $(\Phi_n)_{n\in\N}$  be a subadditive sequence of random continuous functions and let $K$ be a forward invariant random compact set, then 
 \btm
 \item[(1)] the function $\mu\mapsto \lambda(\mu,\varphi)$ is upper semi-continuous from $\Ip$ to $\R$ 
 \item[(2)] the random sequence $(\Phi^K_n)_{n\in\N}$ is subadditive and 
 \item[(3)] if  $\hat{\Phi}^K:= \inf_{n\in\N}\frac{1}{n}\E(\Phi_n^K)= \lim_{n\rightarrow\infty}\frac{1}{n}\E(\Phi^K_n)$ exists, then we have 
 \begin{equation}
 \hat{\Phi}^K=  \sup\{\lambda(\mu,\varphi): \quad \mu\in\mathcal{M}_{\p}(K)\}; 
 \end{equation}
 and the supremum is attained at some $\mu^*\in\mathcal{M}_{\p}(K).$
  \etm
 \end{lem}
 
 \begin{theorem}[Random semiuniform ergodic theorem \cite{Stark13}]\label{semiuni}\rm
 Let $\Theta : \Om\times X\rightarrow \Om\times X$ be a skew product dynamical system with the ergodic base $(\Om, \F,\p, \theta).$ Suppose that $K$ is  a forward invariant random compact set and that $(\Phi_n)_{n\in\N}$ is a subadditive sequence of random continuous function with $\vert \Phi_n\vert^{K}\in L^1(\p)$ for all $n\in\N$ .
 Further, suppose that $\lambda\in\R$ is such that $\lambda(\mu,\varphi)<\lambda$ for all $\mu\in\mathcal{M}_{\p}(K).$  
 
 Then, there exists $\lambda^{\prime}<\lambda$ and an adjusted \footnote{A random variable $C: \Om\rightarrow\R$ is adjusted to $\theta$ if $\lim_{\vert t\vert\rightarrow\infty}\frac{1}{\vert t\vert}C(\theta_t\om) =0,$ or  if and only if the random variable $C\circ\theta-C$ has a $\p$-integrable minorant.} random variable $C:\Om\rightarrow \R$ such that 
 \begin{equation}\label{semi1}
 \Phi_n(\om,x)\leq C(\om)+n\lambda^{\prime}, \quad \text{for all $n$}\in\N, \quad \p-a.e.\om\in\Om\quad \text{and all $x$}\in K(\om).
 \end{equation}
 
 In particular, for $\delta\in(0,\lambda-\lambda^{\prime})$ there is $N(\om)\in\N$ such that 
 \begin{equation}\label{semi2}
 \frac{1}{n}\Phi_n(\om,x)\leq \lambda-\delta, \quad \text{for $n$}\ge N(\om)\quad \text{and $x$}\in K(\om).
 \end{equation}
 \end{theorem}
 
 \begin{theorem}[Variant of random semiuniform ergodic theorem \cite{Stark13}]\label{semiuni2}\rm
 In the situation of Theorem  \ref{semiuni}, there exist $\Lpr<\lambda$ and $k_0\ge 1$ such that for all $k\ge k_0$ there are adjusted random variable $\hat{C}_k:\Om\rightarrow [0,\infty)$ and an ergodic component $\Om_k$ of $\theta^k$ with $\p(\Om_k)\ge \frac{1}{k}$ such that 
 \begin{equation}\label{34}
 \Phi_k(\om,x)\leq \hat{C}_k(\theta^k\om)-\hat{C}_k(\om)+k\Lpr, \quad \text{for $\p$}-a.e.\om\in\Om_k\quad\text{and all $x$} \in K(\om).
 \end{equation}
 The random variables $\hat{C}_k$ can also be choosen to take values in $(-\infty,0],$ furthermore,
 \btm
 \item[(a)] If $(\Phi_n)_{n\in\N}$  is additive, then $k_0=1,$ so that (\ref{34}) holds for $k=1$ and $\p-a.e.$ $\om\in\Om;$
 \item[(b)] If $\theta$ is totally ergodic\footnote{$\theta$ is totally ergodic,  if $\theta^k$ is ergodic for all $k\in\N.$} then (\ref{34}) holds for $\p-a.e.$ $\om\in\Om.$ 
 \etm
 \end{theorem}
 \vspace{.5cm}
 Now let $\varphi$ be a $C^1$ random dynamical system, cocycle property and chain rule of differentiation yield the Jacobian
 \begin{equation}\label{361}
D\varphi(n,\om,x) = \begin{cases}
D\varphi(1,\Theta^{n-1}(\om,x))\cdots D\varphi(1,\om,x)\quad & n\ge 1,\\
x \quad & n=0,\\
(D\varphi(1,\Theta^n(\om,x)))^{-1}\cdots (D\varphi(1,\Theta^{-1}(\om,x)))^{-1}\quad & n\leq -1,
\end{cases}
\end{equation} 
and it follows that $D\varphi$ is a linear cocycle over $\Theta$. \\
Define a random continuous function by $\Phi_n(\om,x):= \log\Vert D\varphi(n,\om,x)\Vert$ and then by (\ref{361}) we have
\begin{align*}
\log\Vert D\varphi(n+m,x)\Vert&= \log\Vert D\varphi(n,\theta^m\om,\varphi(m,\om,x))D\varphi(m,\om,x)\Vert\\
&\leq \log\Vert D\varphi(n,\theta^m\om,\varphi(m,\om,x))\Vert+\log\Vert D\varphi(m,\om,x)\Vert.
\end{align*}
So we have $$\Phi_{n+m}(\om,x)\leq \Phi_m(\om,x)+\Phi_n(\Theta^m(\om,x)),$$
that is to say that the random continuous function $\Phi: \Z\times\Om\times X\rightarrow \R$ is subadditive with respect to $\Theta = (\theta,\varphi).$

We identify the tangent bundle $T\Rd$ of $\Rd$ with $\Rd\times\Rd$ and denote by
\begin{equation}
\lambda(\om,x,v):= \lim_{n\rightarrow\infty}\frac{1}{n}\log\Vert D\varphi(n,\om,x)v\Vert,
\end{equation}
the Lyapunov exponent of $0\neq v\in T_x\Rd\cong\Rd.$
\begin{theorem}[Multiplicative Ergodic Theorem \cite{Arnold1}, \cite{Arnoldan1}, \cite{Arnold}]\rm Supose that
\begin{equation}
\log^{+}\Vert (D\varphi(1,.,.))^{\pm}\Vert\in L^1(\mu).
\end{equation}
Then there exist integers $1\leq r\leq d,$ $1\leq d_1,\cdots,d_r\leq d$ such that $d_1^{+}+\cdots+d_r^{+} = d,$ $d_i^{+} = d_{r+1-i}^{-}$ and $\{\lambda_1^{+},\cdots,\lambda_r^{+}\}\subset\R,$ $\{\lambda_1^{-},\cdots,\lambda_r^{-}\}\subset\R$ such that $\lambda_1^{+}>\cdots>\lambda_r^{+},$ $\lambda_i^{+}=-\lambda^{-}_{r+1-i},$ $1\leq i\leq r,$ a random family of $\mathcal{G}^{+}:=\sigma(\{D\varphi(n,.,.): n\ge 0\})$-measurable orthogonal projectors $(P_i)_{1\leq i\leq r}$ onto the random linear subspace $U_i$ of dimension $d_i^{+},$ a family of $\mathcal{G}^{-}:= \sigma(\{D\varphi(n,.,.): n\leq -1\})$-measurable orthogonal projector $(Q_i)_{1\leq i\leq r}$ onto the random linear subspace $V_i$ of dimension $d_i^{-},$ $1\leq i\leq r,$ and a set $\Gamma\subset \Om\times \R^d$ of full $\mu$-measure which $(\Theta^n)_{n\in\Z}$-invariant such that for any $(\om,x)\in \Gamma$ the following statements hold:
\btm
\item[(i)]\btm
\item[(a)] $\sum_{i=1}^re^{\lambda_i^{+}}P_i(\om,x)= \lim_{n\rightarrow\infty}[D\varphi(n,\om,x)^*D\varphi(n,\om,x)]^{\frac{1}{2n}},$
\item[(b)] if $P(\om,x)_{\ge i} =P_i(\om,x)+\cdots+P_r(\om,x),$ $U_{\ge i}(\om,x)= U_i(\om,x)\oplus\cdots\oplus U_r(\om,x),$ then $$D\varphi(n,\om,x)U_{\ge i}(\om,x) = U_{\ge i}(\Theta^n(\om,x)), \quad n\in\Z, \quad 1\leq i\leq r,$$
\item[(c)] $v\in U_{\ge i}(\om,x)\setminus U_{\ge i+1}(\om,x)$ if and only if $\lambda^{+}_i = \lim_{n\rightarrow\infty}\frac{1}{n}\log \Vert D\varphi(n,\om,x) v\Vert.$
\etm
\item[(ii)]
\btm
\item[(a)] $\sum_{i=1}^re^{\lambda_i^{-}}Q_i(\om,x)= \lim_{n\rightarrow\infty}[D\varphi(-n,\om,x)^*D\varphi(-n,\om,x)]^{\frac{1}{2n}},$
\item[(b)] if $Q(\om,x)_{\ge i} =Q_i(\om,x)+\cdots+Q_r(\om,x),$ $V_{\ge i}(\om,x)= V_i(\om,x)\oplus\cdots\oplus V_r(\om,x),$ then $$D\varphi(-n,\om,x)V_{\ge i}(\om,x) = V_{\ge i}(\Theta^{-n}(\om,x)), \quad n\in\Z, \quad 1\leq i\leq r,$$
\item[(c)] $v\in V_{\ge i}(\om,x)\setminus U_{\ge i+1}(\om,x)$ if and only if $\lambda^{-}_i = \lim_{n\rightarrow\infty}\frac{1}{n}\log \Vert D\varphi(-n,\om,x) v\Vert.$
\etm
\item[(iii)]
\btm
\item[(a)] if $E_i(\om,x) = U_{\ge i}(\om,x)\bigcap V_{\ge r+1-i}(\om,x),$ then $\dim E_i(\om,x) = d_i^{+},$ $\Rd= E_1(\om,x)\oplus\cdots\oplus E_r(\om,x),$ and $D\varphi(n,\om,x)E_i(\om,x) = E_i(\Theta^n(\om,x)),$ $1\leq i\leq r,$
\item[(b)] $0\neq v\in E_i(\om,x) $ if and only if $\lambda_i^{+} = \lim_{n\rightarrow\pm\infty}\frac{1}{n}\log\Vert D\varphi(n,\om,x)v\Vert.$
\etm
\etm
\end{theorem}

The numbers $\lambda_1^{\pm}(\om,x),\cdots,\lambda_r^{\pm}(\om,x)$ are Lyapunov exponents associated with the random measure $\mu$ and random subspaces $(E_i)_{1\leq i\leq r}$ of $\Rd$ are called \textbf{Oseledets spaces}. According to the classical formula due to Furstenberg and Khasminskii (see \cite{Arnold}), the top Lyapunov exponent of a linear cocycle can be expressed as a phase average of a function over a projective space. Arnold and Imkeller \cite{ArnoldIm} derive similar formulas for the remaining Lyapunov exponent of the coycle generated by linear stochastic differential equations. However, the fact that the random subspaces or Oseledets spaces $E_i$ depends on the whole history $\F_{-\infty}^{\infty}$ of Wiener process $W_t$ for the case of stochastic differential equations (with exceptions of $E_1$ resp. $E_r$ which are measurable with respect to $\F_{-\infty}^0$ resp. $\F_{0}^{\infty}$), consitutes basic problems, as one would encounter anticipative stochastic integrand. This problem can be resolved by considering some anticipative problems such as random Fubini theorem for  $\mu_{\om}(dx)d\p(\om)$ and Furstenberg-Khasminskii formula through anticipative calculus (see Arnold and Imkeller \cite{ArnoldIm}, Arnold \cite{Arnoldan}). Precisely, this "random Fubini theorem" is the interchange of the integrals "$\circ dW$" and "$\mu_{.}(dx)$"; the basic idea is to approximate Stratonovich integrals by their Riemann sum. Under the name "substitution formula" this random Fubini theorem has a long history for special case of random Dirac measures $\mu_{.}(dx) = \delta_{y(.)}(dx),$ for a random vector $y$ in $\Rd$ (see Nualart and Pardoux \cite{Nulat}, Millet, Nualart and Sanz \cite{Millt} and Donati-Martin \cite{Donati}). 

Suppose there exists an invariant random compact set $K$ and then by lemma \ref{maxly} there exist a measurable function $b:\Om\rightarrow \Rd$ with $b(\om)\in K(\om)$ such that
  \begin{equation}\label{maxly2}
 \log\Vert D\varphi(n,\om,b(\om))\Vert =\max\{\log\Vert D\varphi(n,\om,x)\Vert: x\in K(\om)\}=:\Phi^K_n(\om).
  \end{equation}
  
Multiplicative ergodic theorem tells us that the numbers $(\om,x)\mapsto \lambda_1^{\pm}(\om,x),\cdots,\lambda_r^{\pm}(\om,x)$ are $\Theta$-invariant and thus constant if $\mu$ is a $\Theta$-ergodic measure.

Put 
\begin{align}
\lambda_i^{+}(\mu,\varphi)&:=\int_{\Om\times\Rd}\lambda_i^{+}(\om,x)d\mu_{\om}(x)d\p(\om),\quad  \mu\in\Ip,\\
\lambda_i^{-}(\mu,\varphi)&:=\int_{\Om\times\Rd}\lambda_i^{-}(\om,x)d\mu_{\om}(x)d\p(\om), \quad  \mu\in\Ip.
\end{align} 
And then the numbers 
 \begin{align}
\lambda_i&:=\sup\{\int_{\Om\times\Rd}\lambda_i^{+}(\om,x)d\mu_{\om}(x)d\p(\om): \mu\in\Ip\},\\
\gamma_i&:=\inf\{\int_{\Om\times\Rd}\lambda_i^{-}(\om,x)d\mu_{\om}(x)d\p(\om): \mu\in\Ip\}
\end{align}
 known as \textbf{extremal exponents} \cite{Hans3} exist at some $\mu^*\in \SM(K),$ which means that $\lambda, \gamma \in\R$.
 We now realise the major hypothesis of the random semiuniform ergodic theorem (Theorem \ref{semiuni}), namely: $\lambda\in\R$ such that 
 \begin{equation}
 \lambda(\mu,\varphi)<\lambda, \quad \text{for all $\mu$}\in \SM(K),
 \end{equation}
 where $\Phi_n(\om,x)= \log\Vert D\varphi(n,\om,x)\Vert
 $ and 
 \begin{align*}
  \lambda(\mu,\varphi)&=\lim_{n\rightarrow\infty}\int_{\Om\times \Rd}\frac{1}{n}\Phi_n(\om,x)d\mu_{\om}(x)d\p(\om)\\
 &=\lim_{n\rightarrow\infty}\int_{\Om\times \Rd}\frac{1}{n}\log\Vert D\varphi(n,\om,x)\Vert d\mu_{\om}(x)d\p(\om).
 \end{align*}
 
\section{Random invariant periodic curves}
We consider a double skew product $H: \T\times\Om\times\Sc\times\R^d\rightarrow\Om\times\Sc\times\Rd,$ \begin{equation}\label{pesk}
H(t,\om,s,x)= (\theta_t\om, s+t\mod 1, \varphi(t,\om,s,x)).
\end{equation}
It is easy to verify that $H$ is a dynamical system on $\Om\times S^1\times \R^d.$ The extra assumption on the random dynamical system $\varphi$ is that $x\mapsto\varphi(t,\om,s,x)$ is at least $C^1$ for all $(\om,s)\in \Om\times S^1.$ 

\begin{definition}[\cite{Zhao}]\rm
Let $\phi(\om):\R\rightarrow\R^d$ be a continuous periodic function with period $\tau\in\N$ for each $\om\in\Om.$ Define $L^{\om}:= \text{graph}(\phi(\om))= \{(s\mod 1, \phi(\om,s)): s\in \R\}.$ If $L^{\om}$  is invariant with respect to the RDS $\varphi$ on $S^1\times \R^d,$ that is, $\varphi(t,\om)L^{\om}= L^{\theta_t\om}$ and there is a minimum $T>0$ (maximum $T<0$) such that for any $s\in[0,\tau)$
\begin{equation}\label{randp}
\varphi(T,\om, s\mod 1, \phi(\om,s)) = (s\mod 1, \phi(\theta_T\om,s))
\end{equation}  for almost all $\om\in\Om,$ then  it is said that $\varphi$ has random periodic solution of period $T$ with random periodic curve $L^{\om}$ of winding number $\tau.$
\end{definition}
 
\subsection*{Assumptions}
Here we give assumptions closely related to that due to Zheng and Zhao \cite{Zhao} on their work on random periodic solutions. 

The average Lyapunov exponent of $\varphi$  associated with $H$-invariant measure $\mu$ is defined by 
\begin{equation}
\lambda(\mu,\varphi)= \lim_{n\rightarrow\infty}\frac{1}{n}\int_{\Om\times\Sc\times\Rd}\Phi_n(\om,s,x)d\mu_{\om}(s,x)d\p(\om),
\end{equation}
where the subadditive random continuous function $\Phi_n$ is given by $$\Phi_n(\om,s,x) = \log\Vert D_x\varphi(n,\om,s,x)\Vert.$$
From the previous section we have seen that for a forward invariant random compact set $K,$ there exists $\mu^*\in\SM(K)$ such that 
$$
\lambda(\mu^*,\varphi) 
=\lim_{n\rightarrow\infty}\frac{1}{n}\int_{\Om\times\Sc\times\R^d}\log\Vert D_x\varphi(n,\om,s,x)\Vert d\mu_{\om}^*(s,x)d\p= \lim_{n \rightarrow\infty}\frac{1}{n}\E\bigg(\max_{(s,x)\in K(\om)}\log\Vert D_x\varphi(n,\om,s,x)\Vert\bigg).
$$
Our assumptions would be based on more general product dynamical system $\Om\times E\ni (\om,s)\mapsto (\theta_t\om,\eta_t(\om,s))$ that is random minimal\footnote{
The product map $\theta\times\eta$ is \textbf{random minimal}, if each forward $\theta\times\eta$-invariant subset $A\subset\Om\times\E$ obeys the following dichotomy:
\btm
\item either $A(\om) = \E$ for $\p-a.e. \om\in\Om$ or
\item  $A(\om)= \emptyset$ for $\p-a.e.\om\in\Om.$
\etm}, where $\E$ is a compact metric space. This condition of minimality is automatic for $\E=\Sc\ni s\mapsto\eta_t(\om,s) = s+t\mod1.$ 

\textbf{Assumption A}
\btm
\item[A1:] Assume that  there exists an invariant random compact set $K$ and that the product metric dynamical system $\theta\times\eta$ is minimal on $\Om\times\E.$
\item[A2:] Assume there is a random invariant compact set $K$ such that $K(\om)$ is connected for $\p-a.e. \om\in\Om.$
\etm 
\textbf{Assumption B}
\btm
\item[B1:] The family $\big((s,x)\mapsto \log\Vert D_x\varphi(k,\om,s,x)\Vert\big)_{\om\in\Om}$ is equicontinuous.
\item[B2:] Given  $\varepsilon>0,$ there exists $r>0$ such that, for all $k\in\T$
\begin{equation}
\sup\{\Phi_k(\om,s,x): \om\in\Om, (s,x)\in B_r(K(\om))\}
\leq \sup\{\Phi_k(\om,s,x): \om\in\Om, (s,x)\in K(\om)\}+\varepsilon.
\end{equation}
\etm
\textbf{Assumption C:} Assume that the Lyapunov exponent $\lambda(\mu,\varphi)$ satisfies $$ \lambda(\mu,\varphi)<0,\quad \text{for all measures $\mu$}\in\mathcal{M}_{\p}(K).$$

Now define $\hq\om:= \theta_{t_1}\om,$ where $t_1$ is the time the particle in $S^1$ rotate a full circle. With this $\hq,$ we reduce our continuous time dynamical system $H$ to discrete time dynamical system $\Hha$. We note that the system $\hq: \Z \times\Om\rightarrow \Om$ is $\p$-preserving such that for any $n,m\in\Z$ one gets
\begin{equation}\label{disch}
H^{n+m}(t_1,\om,s,x)=: \Hha^{n+m}(\om,s,x)= \Hha^n\circ \Hha^m(\om,s,x),
\end{equation}
for $(s,x)\in S^1\times \R^d$ and  $\p$-almost all $\om\in\Om.$

\begin{theorem}[Main result]\label{342}\rm Under assumptions A1 or A2, B1 or B2 and C, there exist  $n(\om)\in\N$ continuous periodic fucntions $\phi_1(\om),\cdots,\phi_{n}(\om)$ with periods $\tau_1(\om),\cdots,\tau_{n(\om)}(\om)\in\N$ respectively such that $$K(\om)=\bigcup_{i=1}^{n(\om)}L^{\om}_i,$$ where $L_i^{\om}=\text{graph}(\phi_i(\om))=\{(s,\phi_i(\om,s)): s\in [0,\tau_i(\om))\}, i=1,\cdots,n(\om).$

Moreover, $\tau_i(\theta_{-t}\om)= 
\tau_i(\om)$ for $i=1,\cdots, n(\om)$ and $n(\theta_{-t}\om)= n(\om),$ $\p$-almost surely.
\end{theorem}
\vspace{1cm}
The prove of this result requires rigorous preparations, we present them in two subsections with the following titles:
\btm
\item Existence of  continuous random invariant graph.
\item Random periodicity of the graph.
\etm

\subsection{Existence of continuous random invariant graph }
Next, for $\varepsilon\in (0,r],$ define  $N_{\varepsilon}(\om,s)$ to be the smallest number of open balls $B_{\varepsilon e^{\hat{C}(\om)}}(x)\subset\R^d$ centred at points $x\in K(\om,s)$ that are required to cover $K(\om,s).$
\begin{lem}[\cite{Stark13}]\label{l3.5}\rm
\btm
\item[(1)]The set $K$ is a random compact set over the base $(\Om\times\E, \F\times\Bb(\E)),$ and $\varphi(k,\om,s,K(\om,s))= K(\theta\times\eta)^k(\om,s)$ for $\p-a.e.\om\in\Om$ and $s\in \E.$ 
\item[(2)] The functions $N_{\varepsilon}$ are $\F\otimes\Bb(\E)$-measurable.
\item[(3)]
For for each $\om\in\Om,$ the function $N_{\varepsilon}(\om,.):\E\rightarrow \N$ is upper semicontinuous. 
\etm
\end{lem}
\textbf{Proof.} \btm\item[(1)] For each $(\om,s)\in \Om\times\E,$ the set $K(\om,s)= \{x\in\R^d: (s,x)\in K(\om)\}$ is compact, as it is a section of a compact set $K(\om).$\\ Denote the metric on $\E$ by $\rho$ and define, for each $n>0$ a metric $d_n$ by $$d_n((s,x),(\hat{s},\hat{x})) = \Vert x-\hat{x}\Vert+n\rho(s,\hat{s}).$$ Each $d_n$ generates a product topology on $\E\times\R^d.$ Further, $\om\mapsto d_n((s,x),K(\om))$ is measurable for all $(s,x)\in\E\times\R^d$ and $s\mapsto d_n((s,x),K(\om))$ is continuous for $\om\in\Om$ and $x\in\R^d.$ Which implies that $(\om,s)\mapsto d_n((s,x),K(\om))$ is measurable.\\
One can check that $d_n((s,x),K(\om))$ is independent of $n$, hence  $(\om,s)\mapsto\sup_{n} d_n((s,x),K(\om)) = d(x, K(\om,s)$ is measurable, so that , $K(\om,s)$ is indeed a random compact set.\\
We know that $x\in K(\om,s)$ if and only if $(s,x)\in K(\om),$\\ so that 
$(\eta(\om,s),\varphi(\om,s,x))\in K(\theta\om)$ if and only if
$\varphi(\om,s,x)\in K((\theta\times\eta)(\om,s))= K(\theta\om,\eta(\om,s)).$
\item[(2)] Now, the set $K$ is a random compact set over the base dynamical system $(\Om\times\E, \F\otimes\Bb(\E)),$ so by random selection theorem, there is a squence of measurable maps $a_k:\Om\times\E\rightarrow\R^d$ such that $K(\om,s) = \overline{\{a_k(\om,s)\}_{k\in\N}}$ for all $(\om,s)\in\Om\times\E.$

Let $\varepsilon>0,$ for $n\in\N$ denote by $\mathcal{L}_n$ the family of subsets of $\N$ with $n$ elements. Then the sets $$V_n:= \bigcup_{\mathcal{L}\in\mathcal{L}_n}\bigcap_{k\in\N}\bigcup_{l\in\mathcal{L}}\big\{ (\om,s)\in\Om\times\E: \Vert a_k(\om,s)-a_l(\om,s)\Vert<\varepsilon\big\}$$ are $\F\otimes\Bb(\E)$-measurable, and $N_{\varepsilon}(\om,s)\leq n$ if and only if $(\om,s)\in V_n.$
\item[(3)] 
It is enough to verify that for $\om\in\Om,$ for each $\alpha\in\R,$  the set $U(\om):=\{s\in\E: N_{\varepsilon}(\om,s)<\alpha\}$ is open subset of $\E.$

We know that $K(\om,s)$ is a compact set, by random selection (Theorem \ref{selection}), there exists $a_k:\Om\times\E\rightarrow\Rd$ such that $K(\om,s)=\text{closure}\{a_k(\om,s): k\in\N\}.$\\ Also, $K(\om)= \{(s,x)\in\E\times\Rd: x\in K(\om,s) \};$ so that for any $y\in K(\om)$ there exists $k\in\N$ such that $\Vert (s,a_k(\om,s))-y\Vert< \frac{\varepsilon}{2}.$\\
Consider the system of subsets of $\E$ $$\mathcal{B}(\om)= \bigcup_{k\in\N}\big\{s\in\E: (s,a_k(\om,s))\in B_{\varepsilon}(y)\big\}.$$  The function $s\mapsto a_k(\om,s)$ is continuous by Proposition \ref{Appa4}, so it follows that the function $(s,a_k(\om,s))$ is also continuous, hence the system of subsets $\mathcal{B}(\om)$ consist of preimage of open sets under continuous function. Therefore, $\mathcal{B}(\om)$ is asystem of open subsets of $\E.$ 

Finally, for each $\alpha\in\R,$ $N_{\varepsilon}(\om,s)<\alpha$ if and only if $s\in\mathcal{B}(\om).$ Thus, $U(\om)$ is open. $ \square$

\etm

\begin{prop}\label{Lypun}\rm
Let $K$ be an invariant random compact set such that Assumptions B and C hold true. Then there are $c>0,$ $\delta>0$ $r>0$ such for all $n\in\T$ and almost all $\om\in\Om$
\begin{equation}\label{lyex1}
\Vert D_x\varphi(k,\om,s,x)\Vert \leq ce^{-\delta k}, \quad \forall (s,x)\in B_r(K(\om)).
\end{equation}

\end{prop}
\textbf{Proof.}
By Theorem \ref{semiuni2} there are $\lambda^{\prime}<0,$ $k\ge 1,$ ergodic component $\Om_k$ of $\theta^k$ and an adjusted random variable $\hat{C}:\Om\rightarrow (-\infty, 0]$ such that 
\begin{equation}\label{ess1}
\Phi_k(\om,s,x)\leq \hat{C}(\theta^k\om)-\hat{C}(\om)+\lambda^{\prime}k, \quad \text{for $\p-a.e.\om\in\Om_k$ and all $(s,x)$}\in K(\om).
\end{equation}
So that, by assumption B2, there are $\lambda^{\prime\prime}<0$ and $r>0$ such that 
\begin{equation}\label{essp}
\log\Vert D_x\varphi(k,\om,s,x)\Vert \leq \hat{C}(\theta^k\om)-\hat{C}(\om)+\lambda^{\prime\prime} k,  \quad \text{for $\p-a.e.\om\in\Om_k$ and all $(s,x)$}\in B_r(K(\om)).
\end{equation}
In particular, for $\delta\in (0, -\lambda^{\prime\prime}),$ there is $N\in\T$ such that 
\begin{equation}
\frac{1}{k}\Phi_k(\om,s,x)\leq -\delta, \quad k\ge N\quad \text{for all $(s,x)$}\in B_r(K(\om)),
\end{equation}
$$\Vert D_x\varphi(k,\om,s,x)\Vert \leq e^{-\delta k}, \quad k\ge N, \quad (s,x)\in B_r(K(\om)).$$

Now, let $\tilde{c}= \max_{k\ge 1}\sup\{e^{\delta k}\Vert D_x\varphi(k,\om,s,x)\Vert\}$ and then $\Vert D_x\varphi(k,\om,s,x)\Vert \leq \tilde{c}e^{-k\delta },$ $k\leq N,$ for all $(s,x)\in B_r(K(\om)).$ So that if we take $c= \max\{\tilde{c}, 1\}$  we have that for all $k\in\T$
$$\Vert D_x\varphi (k,\om,s,x)\Vert \leq c e^{-k\delta },\quad \forall (s,x)\in B_r(K(\om)).$$

\begin{prop}[\cite{Stark13}]\label{34}\rm
Let $H$ be a double skew prodcut and K be a random invariant compact set such that assumptions B1 or B2 and C hold true. Then there are $n\in\N,$ a random variable $c:\Om\rightarrow \Rp$ and a non-empty open forward $\theta\times\eta$-invariant
set $A$  such that for $\p-a.e.\om\in\Om$
\btm
\item $\# K(\om,s) = n$ for all $s\in A(\om),$
\item $\sup\{\#K(\om,s): s\in \E\}<\infty$ and 
\item for all $s\in \E,$ any two different points $y,y^{\prime}\in K(\om,s)$ have distance at least $c(\om).$
\etm
\end{prop}
\textbf{Proof.} There exist $\lambda^{\prime}<0,$  $k\ge 1,$ an ergodic component $\Om_k$ of $\theta^k$ and an adjusted random variable $\hat{C}: \Om\rightarrow (-\infty,0]$ such that 
\begin{equation}\label{35}
\Phi_{k}(\om,s,x)\leq \hat{C}(\theta^k\om)-\hat{C}(\om)+\lambda^{\prime}k, \quad \text{for $\p$-}a.e.\om \quad \text{and all $(s,x)$}\in K(\om).
\end{equation}
Hence in view of assumption B2, there are $r>0$ and $\gamma>0$ such that 
\begin{equation}\label{36}
\log\Vert D_x\varphi^{k}(\om,s,x)\Vert \leq \hat{C}(\theta^k\om)-\hat{C}(\om)- \gamma k, \quad \text{for $\p$-}a.e.\om\in\Om_k\quad \text{and all $(s,x)$}\in B_r(K(\om)).
\end{equation}

Now fix $\varepsilon\in (0,r],$ $\om\in\Om_k$ and $s\in\E,$ and denote $N= N_{\varepsilon}(\om,s).$ There are $x_1,\cdots,x_{N}\in K(\om,s)$ such that $K(\om,s)\subset \bigcup_{i=1}^{N}B_{\varepsilon e^{\hat{C}(\om)}}(x_i).$
As $\varphi^{k}(\om,s,K(\om,s))= K((\theta\times\eta)^k)(\om,s)$ it follows that $$K((\theta\times\eta)^k(\om,s)\subset \bigcup_{i=1}^N\varphi^{k}(B_{\varepsilon e^{\hat{C}(\om)}}(x_i))\subset\bigcup_{i=1}^N B_{e^{-\gamma k}\varepsilon e^{\hat{C}(\theta^{(k)}\om)}}(\varphi^{k}(\om,s,x_i)),$$ with points $\varphi^{k}(\om,s,x_i)\in K((\theta\times\eta)^k(\om,s)).$

Hence, 
\begin{equation}\label{3.6} N_{\varepsilon}((\theta\times\eta)^k(\om,s))\leq N_{\varepsilon e^{-\gamma k}}((\theta\times\eta)^k(\om,s))\leq N_{\varepsilon}(\om,s).
\end{equation}

Consider the restricted system $(\theta\times\eta)^{k}|_{{\Om_k}\times \E}$ and denote the normalized probability measure $\p|_{\Om_k}$ by $\p_k.$

By lemma \ref{l3.5} (2) and (3), there is a subset $\tilde{\Om}_k \subset \Om_k$ of full measure such that the random sets $U_{\varepsilon,\alpha}=\{(\om,s)\in\Om_k\times \E: N_{\varepsilon}(\om,s)<\alpha\}$ are open for all $\alpha\in\R$ and $\varepsilon= e^{-p\gamma r},$ with $p\in\N.$ For measurability purposes we restrict to these countable many values of $\varepsilon$ from now on.
Let $$n_{\varepsilon}(\om)= \min\{\alpha\in\N: U_{\varepsilon,\alpha}(\om)\neq \emptyset\}$$ and the measurability of $n_{\varepsilon}$ follows from lemma \ref{3.2} (4). Due to (\ref{3.6}) we have $$n_{\varepsilon}(\theta^k\om)\leq n_{\varepsilon}(\om)\quad \text{for $\p_k$-}a.e. \om\in\Om_k, $$ and thus, ergodicity of $(\theta^k, \p_k)$ implies that all $n_{\varepsilon}$ are constant $\p_k-a.e.$ 
We denote these constants by $n_{\varepsilon}$ as well. By the first inequality in \ref{3.6} we have $n_{\varepsilon}\leq n_{e^{-\gamma}\varepsilon}.$ But second inequality of (\ref{3.6}) gives $n_{e^{-\gamma}\varepsilon}\leq n_{\varepsilon}$ for all $\varepsilon\in (0,r];$ so all $n_{\varepsilon}$ coincide. Denote their common value by $n$

Using (\ref{3.6}) again, we see that the random open set $U_{r,n}$ is $(\theta\times\eta)^k$-invariant and we have $U_{r,n}= U_{\varepsilon,n}$ for all $\varepsilon.$ Similarly, for each integer $m>n$ the set $U_{r,m}$ is a non-empty $(\theta\times\eta)^k$-invariant open random set and $U_{r,m}= U_{\varepsilon, m}$ for all $\varepsilon.$ 

Next, we show that $K(\om),$ $\om\in \Om_k$ intersects each fibre $\{s\}\times K(\om,s)$ in a finite number of points (this means that the cardinality of $K(\om,s)$ denoted by $\#K(\om,s) = n$). As $K(\om)$ is compact for fixed $\om,$ there exists $(s_1,x_1),\cdots, (s_n,x_n)$ such that \begin{equation}\label{5.11}
K(\om)\subset \bigcup_{i=1}^n B_{r/2e^{\hat{C}(\om)}}(s_i,x_i).
\end{equation}
We will show that for any $s\in U_{r,n}(\om)$ the cardinality of $K(\om,s)$ is at most $n.$ Suppose not, there exists $\hat{s}\in U_{r,n}(\om)$ with $\#K(\om,\hat{s})>n.$ Choose $n+1$ distinct points $x_1,\cdots, x_{n+1}\in K(\om,\hat{s})$ and let $$\min_{i\neq j}\Vert x_i-x_j\Vert .$$ Furthermore, for fixed $k\in\N$ such that $re^{\hat{C}(\om)}\alpha^k< a,$ for some $\alpha\in (0,1)$ and choose, for each $i= 1,\cdots, n+1,$ some $y_i\in \big(\varphi^k(\om,\eta^{-k}(\om,\hat{s}))\big )^{-1}\{x_i\}\in K(\om)$ (such $y_i$ exists due the fact that $K(\om)$ is invariant), due to \ref{5.11}  there exists $l\in \{1,\cdots,n\}$ and $i,j\in \{1,\cdots,n+1\}$ such that  $y_i,y_j\in B_{r/2e^{\hat{C}(\om)}}(x_l).$ Hence the distance between two points is less that $re^{\hat{C}(\om)},$ hence from Propostion (\ref{Lypun}) we conclude that 
\begin{align*}
\Vert x_i-x_j\Vert &= \Vert \varphi^k(\om,\eta^{-k}(\om,\hat{s}),y_i)-\varphi^k(\om,\eta^{-k}(\om,\hat{s}),y_j)\Vert\\ & \leq \alpha^kre^{\hat{C}(\om)}< a, 
\end{align*}
contracting the definition of $a.$
Thus, for $\p_k-a.e.$ $\om\in\Om_k,$ the following holds: 
\btm
\item $\#K(\om,s)= n(\om)$ for all $s\in U_{r,n}(\om),$
\item $\sup\{\#K(\om,s): s\in\E\}<\infty$ and 
\item $d(x,x^{\prime})\ge c(\om):= re^{\hat{C}(\om)}$ for all $s\in\E$ and any two different points $x,x^{\prime}\in K(\om,s).$  
\etm

Let $A_k = U_{r,n} = \{ (\om,s)\in\Om_k\times\E: N_r(\om,s)<n\}.$ If $k=1,$ $A= A_k$ satisfies the assertions of the proposition. Otherwise, we let $$A= \bigcup_{i=0}^{k-1}(\theta\times\eta)^i(A_k),$$ then as $\bigcup_{i=0}^{k-1}\theta^i(\Om_k) = \Om$ up to set of $\p$-measure zero, and the map $\eta^i(\om)$ are homeomorphisms, and $$\varphi^i(\om,s,K(\om,s)) = K((\theta\times\eta)^i(\om,s)),$$ the above conclusions carry over to $\p-a.e.$ $\om\in\Om$ as follows: \\
Observe  that for ergodic component $\Om_k$ of $\theta^k,$ we have that $$\bigcup_{i=0}^{k-1}\theta^i(\Om_k) = \Om, \quad\text{up to set of $\p$-measure zero.} $$
By Theorem \ref{semiuni} and assumption B2 , there exist $r>0,$ $\lambda^{\prime\prime}<0,$ $i\in\{0,1,\cdots,k-1\}$ and an adjusted random variable $C$ such that if we let $\hat{\Om}:= \theta^{i}\Om,$ we have  
\begin{equation}\label{3.10}
\log\Vert D_x\varphi^ k(\theta^{-i}\om,s,x)\Vert\leq C(\theta^{k-i}\om)-C(\theta^{-i}\om)+\lambda^{\prime\prime}k, \quad \text{ for all $k\in\N$, $\p-a.e.\om\in\hat{\Om},$ $(s,x)\in B_r(K(\theta^{-i}\om)).$}
\end{equation}
Following same argument as above, we have that for $\hat{\p}-a.e.$ $\om\in \theta^i\Om_k$ 
\begin{equation}\label{3.11}
\#K(\theta^{-i}\om, \eta^i(\theta^{-i}\om,s)) = n(\theta^{-i}\om) \quad\text{for all $s$}\in A_k(\theta^{-i}\om).
\end{equation}
As $\eta^i(\theta)$ and $\varphi^i(\om)$ are homeomorphisms, hence one-to-one and using the invariance of $K(\om,s)$ with respect to $\varphi^i(\om,s),$ we have that  for $\p_k-a.e.$ $\om\in\Om_k$ 
\begin{equation}
\#\varphi^i(\theta^{-i}\om,s,K(\theta^{-i}\om, \eta^{i}(\theta^{-i}\om,s))) = \#K(\om,s);
\end{equation}
so that for all $\p-a.e.$ $\om\in\Om$ 
\begin{equation}
\#K(\om,s) = n(\om), \quad \text{for all $s$}\in A(\om), 
\end{equation}
where $A= \bigcup_{i=0}^{k-1}(\theta\times\eta)^i(A_k).$
\begin{theorem}[\cite{Stark13}]\label{cong}\rm
Let $H$ be a double skew product and $K$ be a random compact $H$-invariant set such that assumptions A1, B1 or B2 and C hold true. Then there are $n\in\N$ and a random variable $c:\Om\rightarrow \Rp$ such that, for $\p-a.e.\om\in\Om$
\btm
\item $\#K(\om,s) = n$ for all $s\in\E,$
\item the map $s\mapsto K(\om,s)$ from $\E$ to $\mathcal{K}(\R^d)$ is continuous and 
\item for all $s\in\E,$ any two different points $y,y^{\prime}\in K(\om,s)$ have distance at least $c(\om).$
\item If we replace assumption A1 with A2, then for $\p-a.e.\om\in\Om,$ the set $K(\om)$ consists of a single continuous continuous graph. 
\etm
\end{theorem}
\textbf{Proof.} 
The idea here is apply proposition \ref{34}, we note that one the highlights of the proof of the proposition is the construction of a non-empty forward $\theta\times\eta$-invariant random set $A\subset \Om\times\E.$ It follows that $A^{c}$ is a backward $\theta\times\eta$-invariant random set and $A^{c}\neq \E$ for $\p-a.e.$ $\om\in\Om.$ Usingf the assumption that $\theta\times\eta$ is a minimal homeomorphism, we have that $A^{c}= \emptyset$ for $\p-a.e.$ $\om\in\Om$ and hence $A(\om) =\E$ for $\p-a.e.$ $\om\in\Om.$

First and second assertions of proposition \ref{34} together with random compactness of $K,$ implies that $\E\ni s\mapsto\phi_i(\om,s)\in K(\om,s),$ $i=1,2,\cdots, n(\om),$ note that we have applied random selection theorem (Proposition \ref{selection}) here. And then as $K(\om,s)$ is the $s$-section of $K(\om),$ we have that $K(\om)= \text{graph}\phi_i$ $i=1,2,\cdots,n(\om).$

It remains to show that $\E\ni s\mapsto \phi_i(\om,s)$ $i=1,2,\cdots, n(\om)$ are continuous,
Recall that for a compact metric space $\E,$ a function define on $\E$ is continuous on if and only if its graph is compact (see proposition \ref{Appa4} in the appendix).

\textbf{If we replace assumption A1 with A2:} By the second assertion of proposition \ref{34} together with connectedness of $K(\om),$ we have that there is a subset $\hat{\Om}\subset \Om$ of full measure such that $K(\om,s)\subset \R^d$ consists of a single point for all $(\om,s)\in \hat{\Om}\times \E.$ As $K\om)\subset \E\times \R^d$ is compact, then application of proposition \ref{Appa4} shows that $K(\om)$ must be the graph of a continuous function $\phi(\om,.): \E\rightarrow \R^d.$
 As $\{\phi(\om,s)\}= K(\om,s)$ is the only possible selection of $K,$ $\phi$ is measurable (Proposition \ref{selection}). $\square$

\subsection{Random periodicity of the graph}
 We define $(\varepsilon,\delta)$-neighbourhood of $(s,x)\in K(\om)$ and $K(\om)$ as follows

 $$B(s,x,\delta, \varepsilon)=\{(\hat{s},\hat{x}): \vert s-\hat{s}\vert\leq \delta, \Vert \phi(\om,s)-\hat{x}\Vert \leq \varepsilon\},$$ 
  $$\text {and for any $s$}\in \Sc,\quad  B_{s}(\delta,\varepsilon)= \bigcup_{(s,x)\in K(\om)}B(s,x,\delta, \varepsilon),$$   $$\text{so that} \quad B_{\varepsilon}(K(\om))= \bigcup_{s\in S^1}B_{s}(\delta,\varepsilon).$$

\begin{rem}\rm \btm
\item[(1)] In fact, we only need to verify that the continuous functions $\phi(\om)$ from Theorem \ref{cong} is indeed random periodic, given that $\E=\Sc\ni s\mapsto\eta_t(\om,s) = s+t\mod 1$.
\item[(2)] If we suppose that A2 hold, then Theorem \ref{342}  would give that $K(\om)$ consists of a single graph of random periodic function, that is, $n(\om)=1.$ So, it is more interesting to provide the proof with Assumption A1.
\etm
\end{rem}
We start by covering the invariant set $K(\om)$ by boxes within which we have semiuniform contraction (contraction in the direction of $x$ ).
Rather than considering a complete covering of the whole set $K(\om),$ we concetrate on a strip $D_{[s^*-\delta,s^*+\delta]}(\om)=\{\{s\}\times K(\om,s): s\in [s^*-\delta,s^*+\delta]\}$. For any $s\in S^1$ define $$D_s(\om):=\{s\}\times K(\om,s)$$ and for any $(s,x)\in D_s(\om),$ let $U(s,x,\delta,\varepsilon)$ be the interior of $B(s,x,\delta,\varepsilon).$ Then for any $s^*\in S^1,$ $\{U(s^*,x,\delta,\varepsilon): (s^*,x)\in D_{s^*}(\om)\}$ is an open covering of $D_{s^*}(\om).$ By compactness of $D_{s^*}(\om),$ a finite subcover, $U(s^*,x^1_{\om},\delta,\varepsilon),\cdots, U(s^*,x^{p_{\om}}_{\om},\delta,\varepsilon),$ could be found.
Define 
\begin{align*}
U^{\om}(s^*,\delta,\varepsilon)&=\bigcup_{i=1}^{p(\om)}U(s^*,x^i_{\om},\delta,\varepsilon),\\
B^{\om}(s^*,\delta,\varepsilon)&=\bigcup_{i=1}^{p(\om)}B(s^*,x_{\om}^i,\delta,\varepsilon).
\end{align*}
It is easy to see that $B^{\om}(s^*,\delta,\varepsilon)$is the closure of $U^{\om}(s^*,\delta,\varepsilon)$ and that $D_{s^*}(\om)\subset B^{\om}(s^*,\delta,\varepsilon).$

We merge overlap boxes $B(s^*,x^i,\delta,\varepsilon)$ (if any) and work with the connected components of $B^{\om}(s^*,\delta,\varepsilon)$ which we will denote by $B_1^{\om}(s^*,\delta,\varepsilon), B_2^{\om}(s^*,\delta,\varepsilon),\cdots,B_{d(\om)}^{\om}(s^*,\delta,\varepsilon)$ and let the minimal distance between any two of them be $\beta^{\om}>0.$
\begin{prop}\label{4.5}\rm
Under the assumption A1, B1 or B2 and C, $K(\om)$ is a union of a finite number of continuous periodic curves
\end{prop}
\textbf{Proof.} Let $\delta$ be independent of $s\in S^1$ and let $M\in\N$ such that $\frac{1}{M}\leq \delta.$ Define $s_m= \frac{m}{M}, $ $m=1,2,\cdots,M.$ Then $\{(s_{m-1}, s_{m+1}): m=1,2,\cdots,M\}$ ( in which $s_{M+1}= s_1, s_0=s_{M}$) covers $S^1.$ By Theorem \ref{cong}, we know that $D_{[s_{m-1},s_{m+1}]}(\om)$ contains finite number of continuous curves and denote their number by $d(\om,m).$  We claim that $d(\om,m)$ is independent of $m$ and the continuous curves on $D_{[s_{m-1},s_{m+1}]}(\om)$ can be extended to $\Sc.$
\btm
\item\textbf{$d(\om,m)$ is independent of $m$:} The interval $[s_{m}, s_{m+1}]$ is contained in nthe interval $[s_{m-1}, s_{m+1}]$ and thus the strip $D_{[s_m,s_{m+1}]}(\om)$ contains $d(\om, m)$ curves. On the other hand, the interval $[s_m,s_{m+1}]$ is contained in the interval $[s_m,s_{m+2}]$ and the strip $D_{[s_m,s_{m+2}]}(\om)$ contains $d(\om,m+1)$ curves, it follows that the stip $D_{[s_m,s_{m+1}]}(\om)$ contains $d(\om, m+1)$ curves. Hence, $d(\om,m) = d(\om,m+1)$ for all $m\in\N,$ so $d(\om,m)$ must be independent of $m$ and we denote it by $d(\om).$ 
\item \textbf{Extension of curves on $D_{[s_{m-1},s_{m+1}]}(\om)$ to $\Sc$}: Denote the $d(\om)$ curves on each strip $D_{[s_{m-1}, s_{m+1}]}(\om)$ by $\phi_{m,1}(\om),\cdots, \phi_{m,d(\om)}(\om).$ Since $d(\om)$ is independent of $m,$ each $\phi_{m,i}$ can be extended to the whole $\Sc.$ More precisely, let us move to the universal cover $\R$ (as used in the definition of $L^{\om}$) and $s_m=\frac{m}{M}$ for all $m\in\N$ and lift the function $\phi_{m,i}(\om)$ to the cover $\R$ by $\hat{\phi}_{m,i}(\om,s) = \phi_{m\mod M, i}(\om,s\mod 1)$ for $s\in [s_{m-1},s_{m+1}].$ 
 
 Now start by defining $\hat{\phi}_{i}(\om):= \hat{\phi}_{0,i}(\om)$ on $[s_{-1},s_1],$ this agrees with some $\phi_{1,j}(\om)$ for unique $j\in \{1,\cdots,d(\om)\}$ on $[s_0,s_1].$  So defining $\hat{\phi}_{i}(\om) = \phi_{1,j}(\om),$ continue in this way to define $\hat{\phi}_{i}(\om)$ uniquely for $s\ge0.$ Similarly, $\hat{\phi}_{i}$ agrees with some $\phi_{1,j^{\prime}}(\om)$ on $[s_{-1}, s_0],$ we continue to the right to define $\hat{\phi}_{i}(\om)$ for $s\leq 0.$  Hence
 $$D_s(\om) = \{(s\mod 1, \phi_i(\om,s)); i\in\{1,\cdots,d(\om)\}\}$$

\etm

The continuous curves on $D_{[s_{m-1},s_{m+1}]}(\om)$ extend to $\Sc$ and we have the following random return map
$$G^{l}(\theta^{-l}\om): D_s(\theta^{-l}\om)\rightarrow D_s(\om)$$ where $G: \Sc\times \Rd \rightarrow \Sc\times\Rd$ is defined by \begin{equation}\label{rdss}
G^{l}(\om,s,x):= (h^{l}(s),\hat{\varphi}^{l}(\om,s,x))
\end{equation}
where $h(s):= s+t_1\mod1$ and $h^{l}(s) = s+ lt_1\mod1,$ $l\in\N.$
For a fixed $s\in\R,$ the sets $D_s(\theta^{-l}\om)$ and $D_s(\om),$ $l\in\N$ are finite containing $d(\om)$ elements: 
\begin{align*}
D_s(\theta^{-l}\om)&= \{(s\mod 1,\phi_i(\theta^{-l}\om,s) ):  i=1,2,\cdots,d(\theta^{-l}\om)\}, \\ D_s(\om)&= \{(s\mod 1,\phi_i(\om,s)): i=1,2,\cdots,d(\om)\},
\end{align*}
 By the finiteness of $D_s(\om)$ and by continuity of $\phi_i(\om),$ we have
\begin{align*}
\phi_i(\om,s+1)&=\phi_{i_1}(\om,s)\\
\phi_i(\om,s+2)&=\phi_{i_2}(\om,s)\\
&\vdots \\
\phi_i(\om,s+d(\om))&=\phi_{i_{d(\om)}}(\om,s)
\end{align*}
\btm
\item If exact one of $i_1,i_2,\cdots,i_{d(\om)}$ is equal to $i;$ say $i_{\tau_i(\om)}=i.$ Then$$\phi_i(\om,s+\tau_i(\om))=\phi_i(\om,s),$$ for any $s\in \R.$ So $\phi_i(\om)$ is a periodic function of period $\tau_i(\om).$

\item If more than one of $i_1,i_2,\cdots,i_{d(\om)}$ is equal to $i.$ Denote $\tau_i(\om)$ the smallest number $j$ such that $i_j=i.$ Consider, $\tilde{\tau}_i(\om)>\tau_i(\om)$ such that $i_{\tilde{\tau}_i(\om)}=i,$ then \begin{align*}
\phi_i(\om,s+\tau_i(\om))&= \phi_i(\om,s)\\
\phi_i(\om,s+\tilde{\tau}_i(\om))&=\phi_i(\om,s)
\end{align*}
But \begin{align*}
\phi_i(\om,s+\tilde{\tau}_i(\om))&= \phi_i(\om,s+\tilde{\tau}_i(\om)-\tau_i(\om)+\tau_i(\om))\\
&= \phi_i(\om,s+\tilde{\tau}_i(\om)-\tau_i(\om))\\
&=\cdots\\
&=\phi_i(\om,s+\tilde{\tau}_i(\om)-k\tau_i(\om)),
\end{align*}
where $k$ is the smallest integer such that $\tilde{\tau}_i(\om)-(k+1)\tau_i(\om)\leq 0.$ Then by definition of $\tau_i(\om),$
$$\tilde{\tau}_i(\om)-k\tau_i(\om)=\tau_i(\om),$$ so $$\tilde{\tau}_i(\om)=(k+1)\tau_i(\om).$$Theerefore $\phi_i(\om)$ is a periodic function of period $\tau_i(\om).$

\item If none of $i_2,i_2,\cdots,i_{d(\om)}$ is equal to $i.$ In this case, at least two of $i_1,i_2,\cdots,i_{d(\om)}$ must be equal. Say $\tau_2(\om)>\tau_1(\om)$ are such integers such that $i_{\tau_1(\om)}= i_{\tau_2(\om)}$ with smallest difference $\tau_2(\om)-\tau_1(\om).$ Then $$\phi_i(\om,s+\tau_1(\om))=\phi_i(\om,s+\tau_2(\om)).$$
Denote $s+\tau_1(\om)$ by $s,$ then $$\phi_i(\om,s)=\phi_i(\om,s+\tau_2(\om)-\tau_1(\om)), \quad \forall s\in\R.$$
Same as in the second case, we can see for all other possible $\tilde{\tau}_2(\om)$ and $\tilde{\tau}_1(\om),$ $\tilde{\tau}_2(\om)>\tilde{\tau}_1(\om)$ and $i_{\tilde{\tau}_2(\om)}=i_{\tilde{\tau}_1(\om)},$ $\tilde{\tau}_2(\om)-\tilde{\tau}_1(\om)$ must be an integer multiple of $\tau_2(\om)-\tau_1(\om).$ So, $\phi_i(\om)$ is a periodic function of period $\tau_2(\om)-\tau_1(\om).$ $\square$
\etm
Now for any $(s,x)\in B_{\varepsilon}(K(\om)),$ for $l\in\N$ denote 
\begin{align*}
h_1(s)&= h^l(s),\\ \varphi_1(\om,s,x)&=
\hat{\varphi}^l(\om,s,x)= \hat{\varphi}^{l-1}(\om,h^{l-1}(s),\hat{\varphi}^{l-1}(\om,s,x)),\\
H_1(\om,s,x,)&:= (h_1(s), \varphi_1(\om,s,x).)
\end{align*}
We know that there are finite number of continuous periodic functions $\phi_1(\om),\cdots,\phi_{n(\om)}$ on $\R.$ Denote their periods by $\tau_1(\om),\cdots,\tau_{n}(\om)$ respectively. So that $$K(\om)=L_1^{\om}\cup\cdots\cup L^{\om}_{n},$$ where $$L_i^{\om}=\text{graph}(\phi_i(\om))= \{(s\mod 1,\phi_i(\om,s)): s\in [0, \tau_i(\om))\},$$ and from the proof of proposition \ref{4.5}, $\tau_1+\cdots+\tau_n= d(\om)$. But $$H_1(\hat{\theta}^{-l}\om, K(\hat{\theta}^{-l}\om))= K(\om).$$ So \begin{equation}\label{unio}H_1(\hat{\theta}^{-l}\om,L_1^{\hat{\theta}^{-l}\om})\cup\cdots\cup H_1(\hat{\theta}^{-l}\om,L^{\hat{\theta}^{-l}\om}_{n(\hat{\theta}^{-l}\om)})= L_1^{\om}\cup\cdots\cup L^{\om}_{n(\om)}.
\end{equation}

Since $L_i^{\hat{\theta}^{-l}\om}$ is a closed curve and $H_1(\hat{\theta}^{-l}\om)$ is a continuous map , one can get that $H_1(\hat{\theta}^{-l}\om,L_i^{\hat{\theta}^{-l}\om})$ is a closed curve. Moreover, since $H_1$ is a homeomorphism, so 
\begin{equation}\label{inter}
H_1(\hat{\theta}^{-l}\om, L_i^{\hat{\theta}^{-l}\om})\cap H_1(\hat{\theta}^{-l}\om, L_j^{\hat{\theta}^{-l}\om}) = \emptyset, \quad \text{when $i$}\neq j.
\end{equation}

Therefore the left hand side of (\ref{unio}) is indeed a union of $n(\hat{\theta}^{-l}\om)$ distinct closed curves and the right hand side of (\ref{unio}) is a union of $n(\om)$ distinct closed curves. Hence for any $i\in\{1,2,\cdots,n(\hat{\theta}^{-l}\om)\},$ there is a unique $j\in\{1,2,\cdots,n(\om)\}$ such that 
\begin{equation}
H_1(\hat{\theta}^{-l}\om, L_i^{\hat{\theta}^{-l}\om})= L^{\om}_j.
\end{equation}

Similarly, for any $t\in\R,$
\begin{equation}\label{40}
K(\theta_{-t}\om)= L_1^{\theta_{-t}\om}\cup\cdots\cup L_{n(\theta_{-t}\om)}^{\theta_{-t}\om}.
\end{equation}
and
\begin{equation}\label{41}
G(t,\theta_{-t}\om, L_1^{\theta_{-t}\om})\cup\cdots\cup G(t,\theta_{-t}\om, L_{n(\theta_{-t}\om)}^{\theta_{-t}\om})= L_1\cup\cdots\cup L_{n(\om)}^{\om},
\end{equation}
 and here $G:\R\times\Om\times\Sc\times\Rd \rightarrow\Sc\times\Rd$ is defined by $$G(t,\om,s,x) = (s+t\mod1,\varphi_t(\om,s,x)).$$
Without confusion, we can re-order $i's$ and denote the unique $j$ by $i,$ so that for each $\om,$ we have that for any $t\in\R,$
\begin{equation}\label{42}
G(t,\theta_{-t}\om,L_i^{\theta_{-t}\om})= L_i^{\om}
\end{equation}
\begin{lem}\rm
For any $t\in\R,$ $\tau_i(\theta_{-t}\om)= \tau_i(\om)$ for any $i=1,2,\cdots,n(\om).$
\end{lem}
\textbf{Proof.} 
\btm \item Consider when $t=kt_1,$ $k\in\N$ (as indicated earlier in this section, $t_1$ is the time the particle in $\Sc$ rotate a full circle) and note that for any $s\in\{0,1,\cdots,\tau_i^{\theta_{-t}\om}\},$
\begin{equation}\label{43}
\pi_{\Sc}G(t,\theta_{-t}\om, (s,\phi_i(\theta_{-t}\om,s)))= s+k.
\end{equation}
 So for $t=kt_1,$ from (\ref{42}) and (\ref{43}), it turns out that 
\begin{align*}
 &\pi_{\Sc}G(t,\theta_{-t}\om, 0, \phi_i(\om,0)= 0+k,\\  &\pi_{\Sc}G(t,\theta_{-t}\om,\tau_i(\theta_{-t}\om), \phi_i(\theta_{-t}\om,\tau_i(\theta_{-t}\om)))=\tau_i(\theta_{-t}\om)\\ &\text{and} \quad
\pi_{\Sc}G(t,\theta_{-t}\om,\tau_i(\theta_{-t}\om),\phi_i(\theta_{-t}\om,\tau_i(\theta_{-t}\om)))-\pi_{\Sc}G(t,\theta_{-t}\om, 0, \phi_i(\om,0)=\tau_i(\om)
\end{align*}
we have that, $\tau_i(\om)=\tau_i(\theta_{-t}\om).$

\item Consider also, the case when $t\in (kt_1, (k+1)t_1),$ $k\in\N$ and for any $s\in\{0,1,2,\cdots,\tau_i^{\theta_{-t}\om}\},$
$$ \pi_{\Sc} G_t(\theta_{-t}\om,s, \phi_i(\om,s))\in(s+k,s+k+1),$$
since  $\phi_i(\theta_{-t}\om,.)$ is periodic with period $\tau_i(\theta_{-t}\om)$ we have \begin{align*}\pi_{\Sc}G(t,\theta_{-t}\om,\tau_i(\theta_{-t}\om,\phi_i(\theta_{-t}\om))&= \tau_i(\theta_{-t}\om)+k \\&= \tau_i(\theta_{-t}\om)+\pi_{\Sc}G(t,\theta_{-t}\om,0,\phi_i(\theta_{-t}\om,0)).\end{align*}
And for each $(t,\om)\in\R\times\Om,$ using the fact that the curve $L_i^{\om}$ is invariant under the homeomorphism $G(t,\om): S^1\times\R^d\rightarrow S^1\times\R^d,$ to get $$\tau_i(\om) =\pi_{\Sc}G(t,\theta_{-t}\om,\tau_i(\theta_{-t}\om,\phi_i(\theta_{-t}\om))-\pi_{\Sc}G(t,\theta_{-t}\om,0,\phi_i(\theta_{-t}\om,0))= \tau_i(\theta_{-t}\om). \quad \square$$ 
\etm

Finally, let $\phi$ represent any $\phi_i$ and $\tau(\om)$ represent any $\tau_i(\om),$ we already know that $\tau(\theta_{-t}\om)= \tau(\om)$ for any $t\in\R$ and define $\hat{t}= kt_1$. 

Then for any $s\in [0,\tau(\om))$ 
\begin{equation}\label{49}
G(\hat{t},\theta_{-\Ht}\om, s,\phi(\theta_{-\Ht}\om,s)) = (s\mod1, \phi(\om,s)):= \hat{L}^{\om}.
\end{equation}

Therefore from (\ref{49}) and the cocycle property of $\varphi,$ we have that for any $s\in [0,\tau(\om))$
\begin{align}\label{50}
G(\Ht+t,\theta_{-\Ht-t}\om,\hat{L}^{\theta_{-\Ht-t}\om}) &= G(t,\theta_{-t}\om,G(-\Ht,\theta_{-\Ht-t}\om,\hat{L}^{\theta_{-\Ht-t}\om}))\nonumber\\ &= G(t,\theta_{-t}\om,\hat{L}^{\theta_{-t}\om}).
\end{align}

This gives us that for any $s\in [0,\tau(\om))$
\begin{equation}\label{51}
G(t+\Ht,\om,s, \phi(\om,s)) = G(t,\theta_{\Ht}\om,s, \phi(\theta_{\Ht}\om,s)))
\end{equation}
for any $t\leq 0.$ Which implies that $G$ has a periodic curve with period $\Ht$ with winding number $\tau$ and there are $n$ such $\phi.$ That is to say $G$ has $n$ random periodic solutions. $\square$     
\section*{Acknowledgement}
I wish to thank Professor Huaizhong Zhao, my PhD adviser for valuable discussions, suggestions and taking time to read through this work.

\appendix
\section{Appendix}
\begin{theorem}[Continuous image of a compact set ]\label{th1} \rm Let $X$ be a compact topological space and $Y$ be a topological space. Let $f:X\rightarrow Y,$ be continuous. Then, the set $f(X)$ is a compact subset of $Y.$ 
\end{theorem}
\textbf{Proof:} The collection $\{f^{-1}(A); A\in\mathcal{A}\}$ is a covering $X$; these sets are open in $X,$ since $f$ is continuous. Hence, there is a finite subcover $$f^{-1}(A_1),\cdots, f^{-1}(A_n).$$ Then the sets $A_1,\cdots,A_n$ cover $f(X).$ $\square$ 
\begin{lem}\label{th2}\rm
 Every compact subspace of a Hausdorff space is closed
\end{lem}
\textbf{Proof.} Let $A$ be a compact subspace of the Hausdorff space $X.$ We shall prove that $X\smallsetminus A$ is open.
Let $x_0\in X\smallsetminus A.$ We have to show that there is a neighbourhood of $x_0$ that is disjoint from $A.$ 
For each $y\in A,$ let us choose disjoint neighbourhoods $U_{x_0}$ and $V_{y}$ of points $x_0$ and $y$ (using Hausdorff condition).
The collection $\{V_y; y\in A\}$ is a covering of $A$ by open sets open in $X,$ therefore $V_{y_1}, \cdots,V_{y_n}$ cover $A.$ The open set $$V =V_{y_1}\cup\cdots\cup V_{y_n} $$ contains $A$ and it is disjoint from the open set $$U = U_{x_0^1}\cap\cdots \cap U_{x_0^n}.$$ For if $z\in V,$ then $z\in V_{y_i},$ for some $i,$ hence $z\notin U_{x_0^i}.$
We have that $U$ is a neughbourhood of $x_0$ disjoint from $A.$  $\square$
\begin{theorem}\label{th3}\rm
 Let $f:X\rightarrow Y$ be a bijective continuous map. If $X$ is compact and $Y$ is Hausdorff, then $f$ is a homeomorphism.
\end{theorem}
\textbf{Proof.} If $A$ is closed in $X,$ then $A$ is compact so $f(A)$ is comapct.
 Since $Y$ is Hausdorff, then $f(A)$ is closed. $\square$
 
 \begin{prop}\label{Appa4}\rm
Let $X$ and $Y$ be metric spaces. Let $f: X\rightarrow Y$ be any map, define its graph by $$G_{f} = \{(x,y)\in X\times Y: y= f(x)\}\subset X\times Y.$$
\btm
\item[(1)] If $f$ is continuous, the $G_f$ is a closed subset of $X\times Y.$ 
\item[(2)] If $X$ is compact, the map $f$ is continuous if and only if, its graph $G_f$ is compact.
\etm
\end{prop}
\textbf{Proof.}
\btm
\item[(1)] Suppose that $G_f\ni(x_n,y_n)\rightarrow (x,y)$ as $n\rightarrow\infty.$ This implies that $y_n= f(x_n)$ by definition of $G_f$ and $x_n\rightarrow x,$ $f(x_n)\rightarrow y$ as $n\rightarrow \infty.$\\
However, $f$ is assumed to be continuous, so $f(x_n)\rightarrow f(x)$ as $n\rightarrow \infty$ and by uniqueness of limit, we have that $y= f(x)$ and hence $(x,y)\in G_f.$
\item[(2)] The set $X$ is compact, implies for $x_n\in X$ there exists $x_{n_k} \rightarrow x$ as $k\rightarrow\infty.$ \\
Assume that $f$ is continuous, we have that $f(x_{n_k})\rightarrow f(x)$ as $k\rightarrow \infty.$ Let $(x_n,y_n)\in G_f,$  we have that $y_n= f(x_n),$ hence there exists $(x_{n_k},y_{n_k})\rightarrow (x,f(x))\in G_f,$ since $G_f$ is closed.

Assume that $G_f$ is compact. Let $C\subset Y$ be closed and suppose $x_n\in f^{-1}(C)$ and $x_n\rightarrow x\in X$ as $n\rightarrow\infty.$ Then $(x_n,f(x_n))\in G_f$ has a convergent subsequence $(x_{n_k},f(x_{n_k}))$ with limit $(x,y)\in G_f$ (since $G_f$) is closed), so $y= f(x)\in Y$ and since $f(x_{n_k})\rightarrow y$ in $Y,$ it implies that $y= f(x)\in C.$ Thus, $x\in f^{-1}(C)$  which is therefore closed and hence $f$ is continuous.
\etm

\begin{lem} [Facts about random sets \cite{Arnold}, \cite{Stark13}, \cite{Conv}]\label{3.2}\rm
\btm
\item[(1)] Let $A$ be a random set, then $\text{int}(A)$ is a random open set.
\item[(2)] If $\theta\times\eta$ is a measurable homeomorphism and $A$ is a random comapct set, then $(\theta\times\eta)(A)$ is a random compact set with fibres $\eta (\theta^{-1}\om,A(\theta^{-1}\om)).$
\item[(3)] If $(A_n)_{n\in\N}$ are random compact sets, then $\cap_{n\in\N}A_n$ is a random compact set with fibres $\cap_{n\in\N}A_n(\om).$
\item[(4)] If $A$ is random open or closed set, then $\pi_1(A)= \{\om\in\Om: A(\om)\neq \emptyset\}$ is $\F$-measurable.
\etm
\end{lem}
\begin{prop}[Random selection \cite{Arnold},\cite{Conv},\cite{Hans}]\label{selection}\rm
 The set valued map $K$ is a random closed set if and only if there exists a sequence $(k_n)_{n\in\N}$ of measurable maps $k_n: \Om\rightarrow X$ such that $$K(\om)= \text{closure}\{k_n(\om):n\in\N\}\quad \text{for all $\om$}\in\Om.$$

In particular, if $K$ is a random closed set, then there exists a measurable selection, that is, a measurable map $k:\Om\rightarrow X$ such that $k(\om)\in K(\om)$ for all $\om\in\Om.$
\end{prop}

\begin{theorem}[Krylov-Bogolyubov procedure for continuous RDS \cite{Arnold}, \cite{Hans}]\label{Kry}\rm Let $\varphi$ be a continuous RDS on a polish space $X,$ and that $\emptyset\neq\Gamma\subset\pro$ is closed, tight, and convex set of random measures such that $\Theta_t\Gamma\subset\Gamma$ for all $t>0.$ Let $(\nu^N)_{N\in\N}$ be an arbitrary sequence in $\Gamma$ and define the sequence of random measures $(\mu_N)_{N\in\N}$ by the scheme
\begin{equation}\label{kb}
 \mu_N = \begin{cases} \frac{1}{N}\sum_{n=0}^{N-1}\Theta_n\nu^N(.), \quad \T\quad\text{discrete},\\ \frac{1}{N}\int_{0}^{N}\Theta_t\nu^N(.)dt, \quad \T \quad \text{continuous.}
         \end{cases}
\end{equation} 
 (Similarly for $N<0$ if $\T$ is two-sided).  The sequence $(\mu_{N})_{N\in\N}$ has a convergent subsequence and every convergent subsequence of $(\mu_{N})_{N\in\N}$ converges in $\Ip.$
\end{theorem}

 \begin{rem}[Remark on the existence of random invariant measures]\rm 
 
 The conditions of the theorem on the existence of random invariant measures for continuous random dynamical systems (Theorem \ref{Kry}) are statisfied, in particular, for $\Gamma = \Pp$ in the case the state space  $X$ is compact. However, for many interesting and relevant RDS the assumption of a compact state is rather restrictive. In this case, $\Pp$ is not tight, and verification of tightness of a given set $\Gamma$ is not completely trivial. Fortunately, in the paper, we have worked on $\Gamma=\SM(K)$ the set of all random measures supported by an invariant random compact set $K$ which is tight.
 \end{rem}
 \newcommand{\kprime}{k^{\prime}}
 \newcommand{\phih}{\hat{\phi}}

\end{document}